\documentclass[MSNbibl,number,citesort,dvips]{arxbj}
\usepackage{graphicx}
\aid{0}
\volume{19}
\issue{4}
\pubyear{2013}
\firstpage{1294}
\lastpage{1305}
\doi{10.3150/12-BEJSP09} 

\makeatletter
\newcommand{\rrvert}{\vert}
\newcommand{\llvert}{\vert}
\newcommand\bbz{\mathbb{Z}}
\newcommand\bbr{\mathbb{R}}
\newcommand\calc{\mathcal{C}}
\newcommand\calxc{\mathcal{X}_\calc}
\newremark{example}{Example}
\newremark{quest}{Question}
\newcommand{\eqref}[1]{(\ref{#1})}
\makeatother

\begin{document}
\begin{frontmatter}

\title{Some things we've learned\break (about Markov chain Monte Carlo)}
\runtitle{Some things we've learned}

\begin{aug}
\author{\fnms{Persi} \snm{Diaconis}\corref{}\ead[label=e1]{diaconis@math.stanford.edu}}
\runauthor{P. Diaconis} 
\address{390 Serra Mall, Stanford, CA 94305-4065, USA. \printead{e1}}
\end{aug}


%
\begin{abstract}
This paper offers a personal review of some things we've learned about
rates of convergence of Markov chains to their stationary
distributions. The main topic is ways of speeding up diffusive
behavior. It also points to open problems and how much more there is to do.
\end{abstract}

%
\vspace*{-3pt}
\end{frontmatter}

\section{Introduction}\label{sec1}\vspace*{-3pt}

Simulation, especially Markov chain Monte Carlo, is close to putting
elementary probability (Feller Volume I-style) out of business. This
was brought home to me recently in an applied study: Lauren Banklader,
Marc Coram, and I were studying ``smooshing cards,'' a widely used
mixing scheme where a deck of cards is slid around on the table by two
hands. How long should the sliding go on to adequately mix the cards?
To gather data, we mixed 52 cards for a minute and recorded the
resulting permutations 100 times. Why wouldn't these permutations be
random? Our first thoughts suggested various tests: perhaps there would
be too many cards that started adjacent that were still adjacent;
perhaps the cards originally close to the top would stay close to the
top; \ldots. We listed about ten test statistics. To carry out
tests requires knowing the null distributions. I~could see how to
derive approximations using combinatorial probability, for example, for
a permutation $\pi$, consider $T(\pi)=\#\{i:|\pi_i-\pi_{i+1}|=1\}$.
This has an approximate $\operatorname{Poisson}(2)$ distribution with a reasonable error
available using Stein's method \cite{barbour,chatterjee}. For $T(\pi)$
the length of the longest increasing subsequence, some of the deepest
advances in modern probability \cite{baik} allow approximation.

Marc and Lauren looked at me as though I was out of my mind: ``But we
can trivially find null distributions by simulations and know useful
answers in an hour or two that are valid for $n=52$.'' Sigh, of course
they are right, so what's a poor probabilist to do?

One way I have found to go forward has been to study the algorithms
used in simulation. This started with an applied problem: to
investigate the optimal strategy in a card game, a programmer had
generated millions of random permutations (of 52) using 60 random
transpositions. I was sure this was too few (and the simulated results
looked funny).\vadjust{\goodbreak} This suggests the math question, ``how many random
transpositions are needed to mix $n$ cards.'' With Mehrdad Shahshahani
\cite{pd81} we proved that $\frac12 n\log n+cn$ are necessary and
suffice to get $e^{-c}$ close to random. For $n=52$, it takes 400--500.
In retrospect, this is indeed using probability to investigate
properties of an algorithm. I've never worried about finding worthwhile
problems since then.

The literature on careful analysis of Markov chain mixing times is
large. A splendid introduction \cite{levin}, the comprehensive \cite
{aldous}, and the useful articles by Laurent Saloff-Coste \cite
{saloff97,saloff04} give a good picture. There are many other schools
that study these problems. Statistical examples (and theorems) can be
found in \cite{jones,rosenthal}; computer science examples are in \cite
{monte}; statistical physics examples can be accessed via \cite{marti}.
I have written a more comprehensive survey in \cite{pd09}.

The preceding amounts to hundreds of long technical papers. In this
brief survey I~attempt to abstract a bit and ask ``What are some of the
main messages?'' I have tried to focus on applied probability and
statistics problems. Topics covered are
\begin{itemize}
\item Diffusive mixing is slow: Section~\ref{sec2}
\item There are ways of speeding things up (deterministic doubling,
nonreversible chains): Section~\ref{sec3}
\item Some speed-ups don't work (cutting the cards, systematic scans):
Section~\ref{sec4}.
\end{itemize}
Of course, the problems are not all solved and Section~\ref{sec5} gives a list
of open questions I~hope to see answered.

\section{Diffusive mixing}\label{sec2}

\begin{figure}[b]

\includegraphics{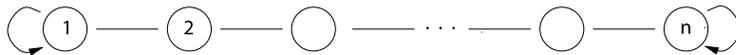}

\caption{Simple random walk on an $n$ point path with $1/2$ holding at
both ends.}
\label{fig1}
\end{figure}

Many Markov chains wander around, doing random walk on a graph. The
simplest example is shown in Figure~\ref{fig1}, a simple random walk on an
$n$-point path.
%

%
\begin{example}
This chain has transition matrix $K(i,j)=1/2,\ |i-j|=1,\
K(1,1)=K(n,n)=1/2$. It has stationary distribution $\pi(i)\equiv1/n$.
Powers of the kernel are denoted $K^l$,
\[
K^2(i,j)=\sum_kK(i,k)K(k,j),\qquad
K^l(i,j)=\sum_kK(i,k)K^{l-1}(k,j).
\]
It is not hard to show that there are universal, positive, explicit
constants $a,b,c$ such that for all~$i,n$,
%
\begin{equation}
ae^{-bl/n^2}\leq\bigl\|K_i^l-\pi\bigr\|\leq
ce^{-bl/n^2} \label{1}
\end{equation}
with $\|K_i^l-\pi\|=\frac12\sum_j|K^l(i,j)-\pi(j)|$.
\label{ex1}
\end{example}\eject

In situations like \eqref{1}, we say order $n^2$ steps are necessary
and sufficient for mixing. The $n^2$ mixing time is familiar from the
central limit theorem which can indeed be harnessed to prove~\eqref{1}.
The random walk wanders around taking order $n^2$ steps to go distance
$n$. This is diffusive behavior.

The same kind of behavior occurs in higher dimensions. Fix a dimension
$d$ and consider the $d$-dimensional lattice $\bbz^d$. Take a convex
set $\calc$ in $\bbr^d$ and look at $\calxc$, the lattice points inside
$\calc$. A~random walk proceeds inside $\calxc$ by picking a nearest
neighbor uniformly at random (probability $1/2d$). If the new point is
inside $\calxc$ the walk moves there. If the new point is outside
$\calxc$ the walk stays. This includes a standard algorithm for
generating a random contingency table with fixed row and column sums:
from a starting table $T$, pick a pair of rows and a pair of columns.
This delineates four entries. Try to change these by adding and
subtracting 1 in pattern $
{{+\atop -}\enskip {-\atop+}}
$ or $
{{-\atop +}\enskip {+\atop-}}
$. This doesn't change the row or column sums. If it results in a table
with nonnegative entries, make the change; otherwise stay at $T$. See
\cite{pd95,pd98} for more on tables.

Returning to the lattice points inside a general convex set, one
expects a bound such as \eqref{1} with $l/n^2$ replaced by $l/(\operatorname
{diam})^2$ for diam the diameter of $\calc$ (length of longest line
inside $\calc$). Theorems like this are proved in \cite{pd11,pd96nash}.
Note that the constants $a,b,c$ depend on the dimension $d$. They can
be as bad as $d^d$, so the results are not useful for high-dimensional
problems. The techniques used are Nash and Sobolev inequalities. There
are extensions of these called \textit{log-Sobolev inequalities} \cite
{pd96sobo,ane} which give good results in high-dimensional problems.
Unfortunately, it is hard to bound the log-Sobolev constant in natural problems.

It is natural to wonder about the choice of the total variation norm $\|
K_i^l-\pi\|$ in \eqref{1}. A variety of other norms are in active use:
\[
\begin{array} {cc@{\quad}l} \bullet&\displaystyle\chi_i^2(l)=\sum
\bigl(K^l(i,j)-\pi(j) \bigr)^2 /\pi
(j)&l^2\mbox{-norm}
\\[6pt]
\bullet&\displaystyle\max_j1-\frac{K^l(i,j)}{\pi(j)}&\mbox{separation}
\\[6pt]
\bullet&\displaystyle\max_j\biggl\llvert 1-\frac{K^l(i,j)}{\pi(j)}\biggr\rrvert
&l^\infty\mbox {-norm}
\\[6pt]
 \bullet&\displaystyle\sum_j\pi(j)\log
\frac{K^l(i,j)}{\pi(j)}&\mbox{Kullback--Liebler.} \end{array} %
\]
One of the things I feel I contributed is this: the choice of distance
doesn't matter; just choose a convenient one and get on with it. Once
you have figured out how to solve the problem with one distance, you
usually have understood it well enough to solve it in others. There are
inequalities that bound one distance in terms of others \cite
{saloff97,gibbs}. The standard choice, total variation, works well with
coupling arguments. Indeed, the maximal coupling theorem says that
there exist coupling times $T$ so that
\[
\bigl\|K_i^l-\pi\bigr\|_{\mathrm{TV}}=P\{T>l\}\qquad\mbox{for
all }l.
\]
The $l^2$ distance works well with eigenvalues. Indeed, for reversible
chains, on a state space of size $n$,
\[
\chi_i^2(l)=\sum_{j=1}^n
\lambda_j^{2l}\psi_j^2(i)
\]
where $\lambda_j,\psi_j$ are the eigenvalues and vectors. Furthermore,
$l^2$ distances allow comparison while total variation doesn't; see
\cite{dyer06,saloff97}. Here ``comparison'' refers to a set of
techniques where a sharp analysis of one chain can be effectively
harnessed to give a useful analysis of a second chain of interest. For
example, on the symmetric group $S_n$, the random transpositions chain
was given a sharp analysis using character theory to show that $\frac12
n\log n$ steps are necessary and sufficient for mixing. From this, the
nonreversible chain ``either switch the top two cards or cut the top
card to the bottom'' was shown to mix in $n^3\log n$ steps. Comparison
uses $l^2$ tools of Dirichlet forms and eigenvalues.

In summary, diffusive behavior occurs for simple random walk Markov
chains on low-dimensional spaces. It leads to unacceptably slow mixing.
The next section suggests some fixes.

\section{Methods of speeding things up}\label{sec3}

The main point made here is that it is often possible to get rid of
diffusive behavior by inserting some simple deterministic steps in the
walk. This is not a well developed area but the preliminary results are
so striking that I hope this will change.
%
\begin{example}[(Uniform distribution on $p$ points)]
Let $p$ be a prime and $C_p$ be the integers modulo~$p$. Simple random
walk goes from $j\in C_p$ to $j\pm1$. It is convenient to change this
to $j\to j,j+1,j-1$ with probability $1/3$. From the arguments in \ref
{sec1} this Markov chain has a uniform stationary distribution $\pi
(j)=1/p$ and from any starting state, order $p^2$ steps are necessary
and sufficient to be close to random. There is diffusive behavior.
\label{ex2}
\end{example}

Consider the following variation: set $X_0=0$ and
\[
X_{n+1}=2X_n+\varepsilon_{n+1} \pmod p
\]
with $\varepsilon_n=0,+1,-1$ with probability $1/3$. This has the same
amount of randomness but intersperses deterministic doubling. Let
$K_n(j)=P\{X_n=j\}$. In \cite{chung87} it is shown that the doubling
gives a remarkable speed-up: order $\log p$ steps are necessary and
sufficient for almost all $p$. One version of the result follows.
\begin{thm*}[(\cite{chung87})]
For any $\varepsilon>0$, and almost all odd $p$, if $l>(C^*+\varepsilon)\log_2p$ then $\|K_l-\pi\|<\varepsilon$ where $C^*= (1-\log_2 (\frac
{5+\sqrt{17}}{9} ) )^{-1}=1.01999186\ldots.$
\end{thm*}

In a series of extensions, Martin Hildebrand \cite
{hildebrand05,hildebrand08} has shown this result is quite robust to
variations: $p$ need not be prime, the probability distribution of
$\varepsilon_i$ can be fairly general, the multiplier 2 can be replaced by
a general $a$ and even $a_{n+1}$ chosen randomly (e.g., 2 or $1/2$ (mod
$p$) with probability $1/2$). The details vary and the arguments require
new ideas.

Once one finds such a phenomenon, it is natural to study things more
carefully. For example, is ``almost all $p$'' needed? In \cite{chung87}
it is shown that the answer is yes: there are infinitely many primes
$p$ such that $\log(p)\log\log(p)$ steps are necessary and sufficient.
Hildebrand \cite{hildebrand09} shows that one cannot replace $C^*$ by 1
in the theorem. In \cite{pd92} similar walks are studied on other groups.

I have heard several stories about how adding a single extra move to a
Markov chain speeded things up dramatically. This seems like an
important area crying out for development. For example, in the
``lattice points inside a convex set $\calxc$'' of Section~\ref{sec1}, is there
an analog of deterministic doubling which speeds up the (diam)$^2$
rate? The reflection walks of \cite{krauth} for the original Metropolis
problem of random placement of non-overlapping hard discs in a box is
an important speed-up of local algorithms. Can it be abstracted?
%
\begin{figure}

\includegraphics{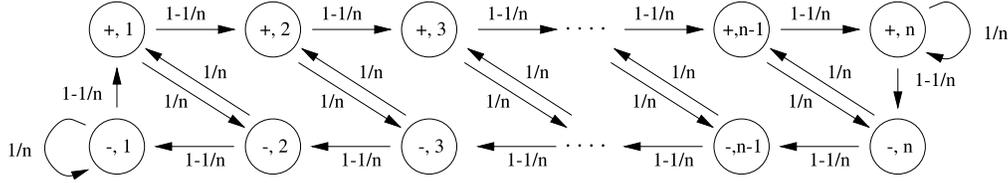}

\caption{A discrete version of hybrid Monte Carlo.}
\label{fig2}
\end{figure}

%
\begin{example}[(Getting rid of reversibility)]
Consider again generating a random point in $\{1,2,3,\ldots,n\}$ by a
local algorithm. In joint work with Holmes and Neale \cite{pdsh00} the
algorithm of Figure~\ref{fig2} was suggested. Along the top, bottom, and side
edges of the graph, the walk moves in the direction shown with
probability $1-(1/n)$. On the diagonal edges the walk moves (in either
direction) with probability $1/n$. The loops indicate holding with
probability $1/n$. While this walk is definitely not reversible, it is
doubly stochastic and so has a uniform stationary distribution.
Intuitively, it moves many steps in one direction before switching
directions (with probability $1/n$). In \cite{pdsh00} it is shown that
this walk takes just $n$ steps to reach stationarity (and this is best
possible for such a local algorithm). The analysis shows that this is a
hidden version of the $X_{n+1}=a_{n+1}X_n+\varepsilon_n$ walk with
$a_{n+1}=1$ or $-1$ with probability $1-(1/n)$ and $1/n$. The walk was
developed as a toy version of the hybrid Monte Carlo algorithm of
lattice field theory~\cite{duane87}. See \cite{neal93} for its
developments in statistics. This is a general and broadly useful class
of algorithms that have resisted analysis. Someone should take up this
challenge!
\label{ex3}
\end{example}

There has been some further development of the ideas in \cite{pdsh00}.
Chen, Lov{\'a}sz and Pak \cite{chen} abstracted the idea to a
``lifting'' of general Markov chains. They showed that the square-root
speed-up (order $n^2$ to order $n$ in the example) was best possible
for their class of algorithms. Hildebrand \cite{hildebrand09} studied
the lifted version of the Metropolis algorithm (based on nearest
neighbor random walk on $\{1,2,\ldots,n\}$) for a general stationary
distribution. The algorithm of Figure~\ref{fig2} chooses to reverse with
probability $1/n$. What about $\theta_n/n$? Evidence in \cite{pdsh00}
suggests that $\theta_n=\sqrt{\log n}$ is better. Gade and Overton \cite
{gade} set this up as an optimization problem, seeking to find the
value of $\theta_n$ that maximizes the spectral gap. In a final
important development, Neal \cite{neal04} has shown that any reversible
Markov chain can be speeded up, at least in terms of spectral gap, by a
suitable nonreversible variant. See \cite{pd12} for further
developments, to spectral analysis for $2^d$-order Markov chains.

In summary, the results of this section show that real speed-ups of
standard algorithms are possible. These results should have practical
consequences: even if it is hard to prove, it is usually easy to find a
few ``big moves'' that preserve the stationary distribution. For a
survey of approaches to designing algorithms that avoid diffusion, see
\cite{andersen}.

\section{Not all speed-ups work}\label{sec4}

One of the joys of proving things is that, sometimes, things that
``everybody knows'' aren't really true. This is illustrated with three
examples: systematic vs. random scans, cutting the cards, and cooking potatoes.
%
\begin{example}[(Systematic vs. random scans)]
Consider applying the Gibbs sampler to a high-dimensional vector, for
example, generating a replication of an Ising model on an $n\times n$
grid. The Gibbs sampler proceeds by updating one coordinate at a time.
Is it better to be systematic, ordering the coordinates and visiting
each in turn, or is choosing a random coordinate (i.i.d. uniform
choices) better? ``Everybody knows'' that systematic scans are better.
Yet, in the only cases where things can be proved, random scan and
systematic scan have the same rates of convergence.
\label{ex4}
\end{example}

Two classes of examples have been studied. Diaconis and Ram \cite{pd00}
studied generation of a random permutation on $n$ letters from the
Mallow's model,
\[
P_\theta(\sigma)=z^{-1}(\theta)\theta^{I(\sigma)},\qquad0<
\theta\leq1,
\]
with $I(\sigma)$ the number of inversions. Here $\sigma=(\sigma
(1),\sigma(2),\ldots,\sigma(n))$ is a permutation of $n$ and $I(\sigma)$
is the number of $i<j$ with $\sigma(i)>\sigma(j)$. This is ``Mallow's
model through Kendal's tau.'' For $0<\theta<1$ fixed, it has $\sigma=$
identity most likely and falls away from this exponentially. The
Metropolis algorithm forms a Markov chain, changing the current $\sigma
$ to $(i,i+1)\sigma$ if this decreases the number of inversions and by
a coin flip with probability $\theta$ if $I((ij)\sigma)>I(\sigma)$;
otherwise the chain stays at $\sigma$. Here, the systematic scan
proposes $(1,2)$, then $(2,3),\ldots,(n-1,n),(n-2,n-1),\ldots,(1,2)$, say.
The random scan chooses $t$ uniformly and independently each time.
Benjamini, Berger, Hoffman and Mossel \cite{benjamini} show that order
$n^2$ random scan steps suffice for random scan. Diaconis and Ram show
that order $n$ systematic scan steps suffice. Since each systematic
scan costs $2n$ steps, the algorithms are comparable. A number of other
scanning strategies and walks on different groups confirm the finding:
being systematic doesn't help to change the order of magnitude of the
number of steps needed for convergence. Two notable features: the
analysis of \cite{pd00} uses Fourier analysis on the Hecke algebra. The
random scan analysis uses deep results from the exclusion process. Both
are fairly difficult. See \cite{bhakta12} for a different approach to proof.

A different set of examples is considered by Dyer--Goldberg--Jerrum
\cite{dyer08}. They studied the standard algorithm for generating a
random proper coloring of a graph with $c$ colors (adjacent vertices
must have different colors). The algorithm picks a vertex and replaces
the color by a randomly chosen color. This step is accepted if the
coloring is proper. How should vertices be chosen to get rapid mixing?
Systematic scan periodically cycles through the vertices in a fixed
order. Random scan chooses vertices uniformly. Intuitively, systematic
seems better. However, their careful mathematical analysis shows the
two approaches have the same convergence rates.

For Glauber dynamics, for Ising and Potts models on graphs, Yuval Peres
(in personal communication) conjectures that random updates are never
faster than systematic scan, and systematic scan can be faster than
random updates by at most a factor of $\log n$ on an $n$-vertex graph.
A~speed-up of $\log n$ is attained at infinite temperature where
systematic scan needs one round of $n$ updates and random scan needs
$n\log n$ updates; see the opening example of \cite{pd00}. Partial
results in the monotone case are in \cite{peres}, Thm.~3.1,~3.2,~3.3.

The results above are tentative because only a few classes of examples
have been studied and the conclusion contradicts common wisdom. It
suggests a research program; a survey of the literature on scanning
strategies is in \cite{pd00}. At least, someone should find one natural
example where systematic scan dominates.
%
\begin{example}[(``Put your faith in Providence but always cut the cards?'')]
Does cutting the cards help mixing? I find it surprising that the
answer is ``Not really and it can even slow things down.'' To say
things carefully, work on $S_n$ the group of all~$n!$ permutations. A
probability on $S_n$ is $Q(\sigma)\geq0,\ \sum_\sigma Q(\sigma)=1$.
Repeated mixing is modeled by convolution,
\[
Q^{*2}(\sigma)=\sum_\eta Q(\eta)Q\bigl(
\sigma\eta^{-1}\bigr),\qquad Q^{*k}(\sigma )=\sum Q(
\eta)Q^{*k-1}\bigl(\sigma\eta^{-1}\bigr).
\]
The uniform distribution is $U(\sigma)=1/n!$. A random cut $C$ puts
mass $1/n$ on each of the $n$-cycles
${1 \atop i}\enskip{2 \atop i+1}\enskip{\ldots \atop \ldots}\enskip{n \atop i-1},\ 1\leq i\leq n$.
It is easy to see, for any of the distances in Section~\ref
{sec1}, $d(C\ast Q,U)\leq d(Q,U)$. So, in this sense, cutting doesn't
hurt (stay tuned!). But does it help? The answer depends on $Q$. For
$Q$ the usual Gilbert--Shannon--Reeds measure for riffle shuffling
$Q^{*k}$ is close to $U$ for $k=\frac32\log_2n+c$ \cite{bayer}. This is
``about 7'' when $n=52$. For general $n$, Fulman \cite{fulman} proves
that applying $C$ after $Q^{*k}$ does not change the $\frac32\log_2n$
rates of convergence.

%
\begin{figure}

\includegraphics{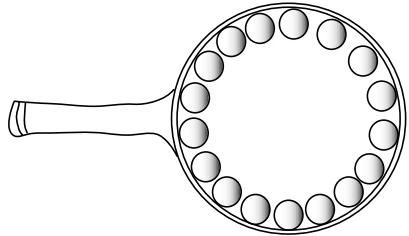}

\caption{16 circular discs inside a pan.}
\label{fig3}
\end{figure}

However, Diaconis and Shahshahani \cite{pd88} construct a probability
measure $Q$ on $S_n$ such that $Q\ast Q=U$ (but $Q\neq U$). For this
$Q$, $(CQ)\ast(CQ)\neq U$. Thus, shuffling twice with this $Q$ gives
perfect mixing but interspersing random cuts fouls things up. Of
course, this $Q$ is not a naturally occurring mixing process. Still, it
shows the need for proof.
\label{ex5}
\end{example}

An example where cutting helps (at least a bit) is in \cite{pd91}.
Here, $Q$ is the random transpositions measure studied by \cite
{pd81,beres,bormashenko}. In \cite{pd81} it is shown that $\frac12n\log
n+cn$ steps are necessary and sufficient for randomness: if $c>0$, $\|
Q^{*k}-U\|\leq2e^{-c}$; if $c<0$, the distance is bounded away from 0
for all $n$. In \cite{pd91}, it is shown that the mixing time of $C\ast
Q$ is $\frac38n\log n+cn$. These are subtle differences. Hard work and
good luck are required to get the lead term constants and cut-off accurately.

%
\begin{example}[(Cooking potatoes)]
When we stir food in a frying pan, e.g., sliced-up potatoes, some
ill-defined ergodic theorem helps to explain why they get (roughly)
evenly browned. One pale mathematical version of this problem considers
$n$ circular discs of potato arranged around the edge of a frying pan
as shown in Figure~\ref{fig3}. Imagine the discs have two sides, heads and
tails. They start with all sides heads-up. At each step, a spatula of
radius $d$ potatoes is inserted at random and all potatoes over the
spatula are turned over in place. For simplicity, assume that $d$ and
$n$ are relatively prime. It is intuitively clear (and not hard to
prove) that with repeated flips, the up/down pattern becomes random;
all $2^n$ patterns are equally likely in the limit.

How long does it take to get close to random, and how does it depend on
$d$? I am surprised that the answer doesn't depend on $d$; a tiny
spatula of diameter 1 or a giant spatula of diameter $n/2$ all require
$\frac14n\log n+cn$ steps (necessary and sufficient) to mix. The result
even holds for ``combs,'' a spatula with teeth that turns over every
other potato among $d$ (or more general patterns).

To see why, regard the potatoes as a binary vector and write $C_2^n$
for the state space. The spatula is a second binary vector, $V$. The
probability measure $Q$ adds a randomly chosen cyclic shift of $V$ to
the current state. Addition is coordinate-wise, mod 2. For
$V=e_1=(1,0,\ldots,0)$, this is just nearest neighbor random walk on the
hypercube, also known as the Ehrenfest urn. The $\frac14n\log n+cn$
answer is well known \cite{pd81}. Consider general $V$. Let $V_1=V,
V_2,\ldots,V_n$ be the $n$-cyclic shifts of $V$. Relatively prime $d$
and $n$ ensures that $V_1,V_2,\ldots,V_n$ form a basis of the space of
binary $n$-tuples. From linear algebra, there is an invertible matrix
$A$ ($n\times n$ mod 2 entries) taking $V_i$ to $e_i,\ 1\leq i\leq n$.
If $0=X_0,X_1,X_2,\ldots$ is the Ehrenfest walk (spatula of size 1) and
$0=Y_0,Y_1,Y_2,\ldots$ is the walk based on $V$, then $P\{Y_k\in S\}=P\{
X_k\in A^{-1}S\}$ for any set $S$. It follow that the total variation
distance to uniformity is the same for the two processes. The same
argument works for any basis $V_1,V_2,\ldots,V_n$ and any distance.
\label{ex6}
\end{example}

Suppose we allow a larger generating set $V_1,V_2,\ldots,V_N$ say with
$N>n$. How should the $\{V_i\}_{i=1}^N$ be chosen to get rapid mixing?
David Wilson \cite{wilson} developed some elegant theory for this question.
\begin{thm*}[(Wilson)]
For all sufficiently large $n$ and $N>n$, and $V_1,V_2,\ldots,V_N\in
C_2^n$, the random walk based on repeatedly adding a uniformly chosen
$V_i$ satisfies
\begin{enumerate}
\item for any choice of $V_1,\ldots,V_N$, if $k<(1-\varepsilon)T(n,N)$ then
$\|Q^{*k}-U\|>1-\varepsilon$;
\item for almost all choices of $V_1,\ldots,V_N$, if $k>(1+\varepsilon
)T(N)$ then $\|Q^{*k}-U\|<\varepsilon$ provided the Markov chain is ergodic.
\end{enumerate}
\end{thm*}

Here $T(n,N)=\frac{N}2\frac1{1-H^{-1}(n/N)}$, $H(x)=x\log_2\frac
1{x}+(1-x)\log_2\frac1{1-x}$, $0\leq x<1$. Note that almost all choices
in item 2 of the theorem will be ergodic when $N-n$ is sufficiently
large. For example, when $N=2n$, $T(n,N)\doteq0.24853n$ steps are
required. Further details are in \cite{wilson}.

\section{Open questions}\label{sec5}

\begin{quest}
In item 2 of Wilson's theorem (Example~\ref{ex6}), the result holds for
almost all choices $V_1,V_2,\ldots,V_N$. Can an explicit set be found,
e.g., for $N=2n$?
\label{quest1}
\end{quest}
%
\begin{quest}
The same set of problems can be considered for any group $G$. If a
generating set $S$ is chosen at random, what is the typical rate of
convergence? This is the topic of random random walks. Hildebrand \cite
{hildebrand05} gives a survey. Babai, Beals, and Seress \cite{babai}
give the best bounds on the diameter of such random Cayley graphs.
These may be turned into (perhaps crude) rates of convergence via
bounds in \cite{pd93}. I cannot resist adding mention of one of my old
conjectures. For the alternating group $A_n$, it is known that a
randomly chosen pair of elements generate $A_n$ with probability
approaching 1. I conjecture that the random walk based on any
generating pair gets random in at most $n^3\log n$ steps.
\label{quest2}
\end{quest}
%
\begin{quest}
Fix a generating set $S\subseteq G$. What element should be added to
$S$ to best speed up mixing? For example, suppose $G=S_n$ (for some odd
$n$) and $S=\{(1,2),(1,2,3,\ldots,n)\}$, a~transposition and an
$n$-cycle. It is known that order $n^3\log n$ steps are necessary and
suffice for randomness \cite{pd93,wilson}. Is there a choice of $\sigma
$ to be added that appreciably speeds this up? For $S_n$, it is
conjectured that all such walks have a sharp cutoff \cite{chen}.
\label{quest3}
\end{quest}
%
\begin{quest}
One may ask a similar question for random walk on any graph. To be
specific, consider a connected $d$ regular graph with $n$ even. Thus,
nearest neighbor random walk has a uniform stationary distribution. Add
in $n/2$ edges forming a perfect matching. This gives a $(d+1)$ regular
graph. What choice of edges give fastest mixing? If the original graph
is an $n$-cycle and thus 2-regular, \cite{chung89} shows that a random
matching improves the diameter to $\log_2n+o(1)$. She gives an explicit
construction of a matching that has diameter $2\log_2n+o(1)$. These
diameter bounds translate into eigenvalue bounds and so bounds on rates
of convergence using standard tools. However, something is lost in
these translations and it would be worthwhile to know accurate rates of
convergence to the uniform distribution.
\label{quest4}
\end{quest}

An important variation: consider a reversible Markov chain $K(x,y)$ on
a finite set $\mathcal{X}$ with stationary distribution $\pi(x)$.
Suppose a weighted edge is to be added to the underlying graph and the
resulting Markov chain is ``Metropolized'' so that it still has
stationary distribution $\pi(x)$. What edges best improve mixing, or
best improve the spectral gap? These questions are closely related to
Section~\ref{sec3}.


%


\printhistory

\end{document}